\documentclass[a4paper]{amsart}
\usepackage{amssymb}
\usepackage[all]{xy}
\usepackage{ifthen}
\newtheorem{theorem}{Theorem}[section]

\newtheorem{lemma}[theorem]{Lemma} 
\newtheorem{fact}{Fact}[section]

\newtheorem{corollary}[theorem]{Corollary}
\theoremstyle{definition}

\newtheorem{definition}{Definition}[section]
\theoremstyle{remark}
\newtheorem{notation}[definition]{Notation}

\newcommand{\field}[1]{\mathbb{#1}}
\newcommand{\bbbn}{\field{N}}
\newcommand{\bbbr}{\field{R}}

\newcommand{\vrightarrow}[1]{\
  \xy\xymatrix{*!{}\ar@{.>}[r]^{\scriptscriptstyle #1}&*!{}}\endxy\ }
\newcommand{\vlrightarrow}[2]{\
  \xy\xymatrix@+#1{*!{}\ar@{.>}[r]^{\scriptscriptstyle #2}&*!{}}\endxy\ }

\DeclareMathOperator\depth{{\rm td}}
\DeclareMathOperator\clos{{\rm clos}}
\DeclareMathOperator\height{{\rm height}}
\DeclareMathOperator\tw{{\rm tw}}

\DeclareMathOperator\mad{{\rm mad}}
\newcommand\vcol[1]{\gamma_{#1}}

\newcommand\ecol[1]{\Upsilon_{#1}}
\newlength\cplxlen
\newlength\cplxlenb
\newcommand{\gdens}[2]{
\settowidth{\cplxlen}{$\scriptstyle #1$}
\settowidth{\cplxlenb}{$\scriptstyle #2$}
\setlength{\cplxlen}{0.5\cplxlen}
\hspace{\cplxlen}\overset{\hspace{-\cplxlen}{#1}\hspace{-\cplxlen}}{\nabla}_{#2}
\addtolength{\cplxlen}{-\cplxlenb}
\ifthenelse{\lengthtest{\cplxlen > 0pt}}{\hspace{\cplxlen}}{}
}
\newcommand{\rdens}[1]{\nabla_{#1}}
\newcommand{\dens}{\nabla}
\newcommand{\card}[1]{\lvert{#1}\rvert}

\DeclareMathOperator\md{\Delta^{\rm --}}

\def\vG{\vec{G}}
\def\vH{\vec{H}}

\def\vP{\vec{P}}
\def\vY{\vec{Y}}
\def\vL{\vec{\Lambda}}
\newcommand{\rdpath}[1]{\xy\xymatrix{*!{}\ar@{~2>}[r]^{#1}&*!{}}\endxy}
\newcommand{\ldpath}[1]{\xy\xymatrix{*!{}&*!{}\ar@{~2>}[l]_{#1}}\endxy}
\newcommand{\dpath}[3]{{#1}\rdpath{#2}{#3}}
\newcommand{\lambd}[5]{{#1}\rdpath{#2}{#3}\ldpath{#4}{#5}}

\begin{document}

\title{Grad  and Classes with Bounded Expansion I. Decompositions.}
\author{Jaroslav Ne\v
  set\v ril}
\thanks{Supported by grant 1M0021620808 of the Czech Ministry of Education}
\address{Department of Applied Mathematics\\
  and\\
  Institute of Theoretical Computer Science (ITI)\\
  Charles University\\
  Malostransk\' e n\' am.25, 11800 Praha 1\\
  Czech Republic} \email{nesetril@kam.ms.mff.cuni.cz} \author{Patrice
  Ossona de Mendez}
\address{Centre d'Analyse et de Math\'ematiques Sociales\\
  CNRS, UMR 8557\\
  54 Bd Raspail, 75006 Paris\\
  France} \email{pom@ehess.fr}
\begin{abstract}
We introduce classes of graphs with {\em bounded expansion} as a
generalization of both proper minor closed classes and degree bounded
classes. Such classes are based on a new invariant, the
{\em greatest reduced average density (grad) of $G$ with rank $r$},
$\rdens{r}(G)$. For these classes we prove the existence of several
partition results such as the existence of low tree-width and low
tree-depth colorings. This generalizes and simplifies several earlier
results (obtained for minor closed classes). 
\end{abstract}

\maketitle
\section{Introduction}
\label{sec:intro}
Let us start with the following particular case which illustrates some
of the motivation of this paper:
It is well known that not only the chromatic number of
planar graphs is bounded
but so are various of its variants such as acyclic or star chromatic number
(by $5$ and $20$, see for instance \cite{AMS} and \cite{albertson04}). For
which other classes of graphs does this hold?  While these variants 
of chromatic number are unbounded even 
for bipartite graphs, we proved
in \cite{Taxi_jcolor} that any proper minor closed class of graphs
has a bounded star chromatic number:
For any minor closed class of graphs $\mathcal C$ excluding at least one
graph --- what we shall call a proper minor closed class ---
there exists an integer $N(\mathcal C)$ such that any graph
$G\in\mathcal C$ has a colorex by $N(\mathcal C$ colors such that any two colors
induce a star forest. Thus also acyclic chromatic number of graphs from a 
proper minor closed class is bounded. This particular case  also 
follows from a 
 recent result of DeVos et al.
 \cite{2tw} who  proved, using the
Structural Theorem of Robertson and Seymour \cite{GM16}, that for any fixed
integer $p\geq 1$, any proper
minor closed class of graphs has a bounded coloring such that any
$i\leq p$ parts induce a graph of tree-width at most $(i-1)$. 
Such a coloring is called low tree-width coloring.

In \cite{Taxi_tdepth}, we presented a strengthened version of
\cite{2tw} : 
we introduced the
 tree-depth
of a graph and proved that
for any fixed $p$, any proper
minor closed class of graphs has a bounded coloring such that any
$i\leq p$ parts induce a graph of tree-depth at most $i$. We also proved
that tree-depth is the best graph invariant with this property 
(see \cite{Taxi_tdepth} and below for more details) .
Also this result uses \cite{2tw} and thus also the 
Structural Theorem. Such a coloring is called low tree-depth coloring
and 
this naturally leads  to  a sequence $\chi_1, \chi_2, \ldots$
of chromatic numbers $\chi_p$,
where $\chi_1$ is the usual chromatic number, $\chi_2$ is the star
chromatic number and, more generally, $\chi_p$ is the minimum number
of colors such that any $i\leq p$ parts induce a graph with tree-depth
at most $i$.

It is well known that $\chi_1$ is bounded on a class of graphs if
the maximum average degree of graphs in the class is
bounded. In \cite{Taxi_jcolor}, we actually proved that $\chi_2$ is bounded if
the graphs obtained by contracting star forests have bounded maximum
average degree. Also, if $\chi_2$ is bounded then so is the maximum
average degree
 (Assume $\chi_2(G) \leq N$. Then for any two colors 
$i\neq j, i, j \leq N$,
orient the edges of $G$  such that any vertex has 
indegree at most one
in the star forest induced by colors $i$ and $j$. Then the indegree of
any vertex 
is at most $\binom{N}{2}$ and thus the graph has maximum average
degree at most $2\binom{N}{2}$.)

This indicates that  the minor closed classes are perhaps  not the most natural 
restriction in the context of graph partitions. One is naturally led to the 
study of minors with bounded depth (of the contracted forest) and their 
edge densities. This in turn leads to the notion of bounded expansion 
which is the central notion of this paper. 

Very schematically
this relationship between the $\chi_p$'s and the bounded depth minors
naturally leads to the following two  questions:

Do there exist integral functions $f_1$ and $f_2$ such that, for any integer $p$:

\begin{itemize}
\item If the minors of depth at most $f_1(p)$ of the graphs of a class
$\mathcal C$ have
bounded maximum average degree then the graphs in $\mathcal C$ have
bounded $\chi_p$,
\item If the graphs in $\mathcal C$ have bounded $\chi_{f_2(p)}$ then all  the
minors of depth at most $p$ of the graphs of a class
$\mathcal C$ have
bounded maximum average degree.
\end{itemize}

In this paper, we prove that  both questions have a  
positive answer. This is the main result of this paper formulated below 
as Theorem \ref{thm:main}. 
It implies the above result of \cite{Taxi_tdepth}. 
Perhaps more interestingly our proof does not rely on the Structural Theorem and yield an effective
algorithm (in fact a linear algorithm, see our companion paper \cite{POMNII}).

Let us describe this development in a greater detail:
The concept of tree-width \cite{Halin},\cite{GM1},\cite{Wagner} is
central to the analysis of graphs with forbidden minors done by
Robertson and Seymour and gained much algorithmic attention
thanks to the general complexity result of Courcelle about
monadic second-order logic graph properties decidability for graphs
with bounded tree-width \cite{Courcelle1},\cite{Courcelle2}.
This computational property (and similar algorithmic aspects), as well
as a question of R. Thomas \cite{2twconj},
motivated the study of graph partitions where $k$ parts induce a
subgraph of tree-width at most $(k-1)$. Such partitions have been proved
to exists by DeVos et al. for proper minor closed classes of graphs
\cite{2tw}, relying on Structural Theorem of Robertson and Seymour on
the structure of graphs without a particular graph as a minor \cite{GM16}. 
This result has been extended by the authors to tree-depth
decompositions in \cite{Taxi_tdepth}. Advancing the definition of 
tree depth let us recall the definition of the tree width by means of 
$k$-trees: A {\em $k$-tree} is a graph which is either a clique of size at most
$k$ or is obtained from a smaller $k$-tree by adding a vertex adjacent
to at most $k$ vertices which are pairwise adjacent. The {\em
tree-width} $\tw(G)$ of a graph $G$ is the smallest integer $k$ such that $G$
is a subgraph of a $k$-tree, that is: a {\em partial $k$-tree}.
The {\em tree-depth} $\depth(G)$ of a connected  graph $G$ is the minimum height
of a rooted tree which closure contains $G$ as a subgraph (height is
defined here as the maximum number of vertices of a path from the root
to a leaf of the tree; the closure of a rooted tree is the graph formed 
by the ancestor relation). (The tree depth of a disconnected graph $G$ is 
the maximal
tree depth of a component of $G$.)

The tree depth is a  minor monotone invariant. It  is related to
the tree-width by $\tw(G)+1\leq\depth(G)\leq\tw(G)\log_2 n$, where $n$ is
the order of $G$ and is actually equal to the {\em vertex ranking
number} \cite{vertex_ranking}\cite{Schaffer} and to the minimum height
of an {\em elimination tree}  \cite{vertex_ranking}. For our purposes it is 
important that $\depth(G)$ has an  alternative
definition by means of  {\em centered coloring}: a  coloring of the vertices
of a graph $G$ is called centered if in any connected subgraph $G'$ of
$G$ some color appears exactly
once (thus a centered coloring is necessarily proper).
It may be seen then that the tree-depth of a graph $G$ is the minimum 
number of colors in a
centered coloring of $G$. As well as graphs with large tree-width 
may be characterized by large  grid minors,
tree-depth may be characterized by excluded paths: a graph has large
tree-depth if and only if it includes a long path. 

Generalizing \cite{2tw} we proved in \cite{Taxi_tdepth} the following:

\begin{theorem}[Corollary 5.3 of \cite{Taxi_tdepth}]
\label{mainold}
For any proper minor closed class of graphs $\mathcal K$ and for any
fixed integer $p\geq 1$, $\chi_p(G)$ is bounded on $\mathcal K$.
\end{theorem}

An alternative way to look at this result is the following: for any
integer $k$ and any proper minor closed class of graphs $\mathcal K$,
there exists an integer $N(\mathcal K,k)$ such that any subgraph
$H\subseteq G$ gets at least $\min(k,\depth(H))$ colors (hence $i<k$ parts
induce graphs of tree-depth at most $i$).

In\cite{Taxi_tdepth} we proved that  this  statement is optimal
in the following sense: Let $\phi$ be an integral graph function
(i.e. we assume that  $\phi(G)$ is an integer  for any graph
$G$). Assume that for  any integer $k$ and for
any proper minor
closed class $\mathcal K$ there exists an integer $N({\mathcal K},k)$
such that any graph $G \in \mathcal K$ has a partition into $\leq
N({\mathcal K},k)$ parts with the property that any subgraph
$H\subseteq G$
gets at at least $\min(k,\phi(H))$ colors. Then $\phi(H) 
\leq \depth(H)$. 

Here we
extend 
Theorem \ref{mainold} to more general classes of graphs. In fact it appears that proper minor closed classes are unnecessary restrictive for the validity of Theorem \ref{mainold}. 

Let $f$ be a function assigning to every positive integer $n$ a real value $f(n)$.
Instead of dealing with proper minor closed classes we shall work with classes of graphs with 
$f$-bounded expansion. This definition is introduced in Section \ref{sec:grad}. Informally, a graph $G$ is 
said to have $f$-bounded expansion if every minor $G'$ of $G$ which we obtain by 
contracting a disjoint union of connected subgraphs of radius $\leq r$
and then deleting some
vertices have edge density bounded by $f(r)$.
The main consequence of our approach here is a generalization of 
Theorem \ref{mainold}
to the classes of graphs with $f$-bounded expansion. This is indeed a generalization 
as each proper minor closed class has expansion bounded by a constant. Also
bounded degree graphs are fitting into this scheme (they are bounded by an exponential function).
(See Section \ref{sec:grad} where  the bounded expansion is defined and discussed
in detail.) Actually, we not only extend Theorem \ref{mainold} to
classes with bounded expansion but prove that it cannot be extended
further: classes with bounded expansion may be actually characterized
by the validity of Theorem \ref{mainold}.
 
It is perhaps surprising that one can prove the full analogy of Theorem \ref{mainold} 
on this level of generality.
The main reason for this is that we approach the decomposition theorem 
via graph orientations and their local properties. Note that 
triangulated graphs, like $k$-trees, have orientations with strong
local properties.
A digraph $\vG$ is {\em fraternally oriented} if $(x,z)\in E(\vG)$ and
$(y,z)\in E(\vG)$ implies $(x,y)\in E(\vG)$ or $(y,x)\in E(\vG)$. This
concept was introduced by Skrien \cite{skrien82} and a
characterization of fraternally oriented digraphs having no symmetrical
arcs has been obtained by Gavril and Urrutia \cite{gavril92}, who also proved that
triangulated graphs and circular arc graphs are all fraternally
orientable graphs. An orientation is {\em transitive} if $(x,y)\in
E(\vG)$ and $(y,z)\in E(\vG)$ implies $(x,z)\in E(\vG)$. It is obvious
that a graph  has an acyclic transitive
fraternal orientation in which every vertex has indegree at most
$(k-1)$ if and only if it is the closure of a rooted forest of height
$k$. It follows that tree-depth and transitive fraternal orientation
are closely related.

This paper is organized as follows: In Sections
 \ref{sec:tw},\ref{sec:td},\ref{sec:grad}
 we introduce the above notions in a greater detail.
The key notion is the notion of the {\em greatest reduced average density (grad)} $\rdens{r}(G)$
of rank $r$ of a graph $G$. We then derive several results about local properties of orientations.
 This is the
reason why we use or introduce relaxed versions, like {\em $p$-centered
colorings} (in which in every subgraph, either some color appears
exactly once or at least $p$ colors appear), or {\em transitive fraternal
augmentations} (each augmentation step consists in adding the missing
arcs while applying the fraternity and transitivity rules on the
initial arcs). 
The Section \ref{sec:lex} is devoted to the proof of the stability of the notion
of classes with bounded expansion with respect to the lexicographic
product with an arbitrary fixed size complete graph (Lemma
\ref{lem:gdens}). This key lemma will allow to prove in Section \ref{sec:aug} the
existence of transitive fraternal augmentations with indegrees bounded
as a function of the grad. These augmentations will be used in Section
\ref{sec:back} to exhibit $p$-centered colorings, eventually leading us to
Theorem \ref{thm:main} in Section \ref{sec:cncl}.

Further corollaries and applications of our method will appear in the
3 companion papers,  
see \cite{POMNII,POMNIII,POMNIV}.

\section{Low tree-width coloring}
\label{sec:tw}

A {\em $k$-tree} is recursively defined as a single vertex graph or a
graph obtained from a smaller $k$-tree by adding a vertex adjacent to a
clique of size at most $k$. The {\em tree-width} $\tw(G)$ of a graph
$G$ is the minimum integer $k$ such that $G$ is a subgraph of a
$k$-tree.

A class $\mathcal C$ has a {\em low tree-width coloring} if, for any
integer $p\geq 1$, there exists an integer $N(p)$ such that any graph
$G\in\mathcal C$ may be vertex-colored using $N(p)$ colors so that
each of the connected components of the subgraph induced
by any $i\leq p$ parts has tree-width at most $(i-1)$.
According to this definition, the result of DeVos et al. may be
expressed as

\begin{theorem}[\cite{2tw}]
\label{th:2tw}
Any minor closed class of graphs excluding at least one graph has a
low tree-width coloring.
\end{theorem}
\section{Low tree-depth coloring and $p$-centered colorings}
\label{sec:td}
In \cite{Taxi_tdepth}, we introduced
the {\em tree-depth} $\depth(G)$ of a graph $G$ as follows:

A {\em rooted forest} is a disjoint union of rooted trees.
The {\em height} of a vertex $x$ in a rooted
forest $F$ is the number of vertices of a
path from the root (of the tree to which $x$ belongs to) to $x$ and is noted $\height(x,F)$.
The {\em height} of $F$ is the maximum height of the vertices of $F$.
Let $x,y$ be vertices of $F$. The vertex $x$ is an {\em ancestor} of $y$ in $F$ if $x$ belongs to
the path linking $y$ and the root of the tree of $F$ to which $y$ belongs to.
The {\em closure} $\clos(F)$ of a rooted forest $F$ is the graph with
vertex set $V(F)$ and edge set $\{\{x,y\}: x\text{ is an ancestor of
  }y\text{ in }F, x \neq y\}$. A rooted forest $F$ defines a partial order on its set of vertices:
$x\leq_F y$ if $x$ is an ancestor of $y$ in $F$.
The comparability graph of this partial order is obviously $\clos(F)$.
The {\em tree-depth} $\depth(G)$ of a graph $G$ is the minimum height
of a rooted forest $F$ such that $G\subseteq\clos(F)$. As a
consequence, we have:

\begin{lemma}[\cite{Taxi_tdepth}]
\label{lem:recur}
Let $G$ be a graph and let $G_1,\dots,G_p$ be its connected components. Then:
\begin{equation*}
\depth(G)=\begin{cases}
1,&\text{if }|V(G)|=1;\\
1+\min_{v\in V(G)}\depth(G-v),&\text{if }p=1\text{ and }|V(G)|>1;\\
\max_{i=1}^p \depth(G_i),&\text{otherwise.}
\end{cases}
\end{equation*}
\end{lemma}

As we introduced low tree-width coloring, we say that
a class $\mathcal C$ has a {\em low tree-depth coloring} if, for any
integer $p\geq 1$, there exists an integer $N(p)$ such that any graph
$G\in\mathcal C$ may be vertex-colored using $N(p)$ colors so that
each of the connected components of the subgraph induced
by any $i\leq p$ parts has tree-depth at most $i$.
As $\depth(G)\geq\tw(G)-1$, a class having a low-tree depth coloring
has a low tree-width coloring.
In \cite{Taxi_tdepth} is proved a strengthening of
Theorem \ref{th:2tw}:
\begin{theorem}[\cite{Taxi_tdepth}]
\label{th:ltd}
Any minor closed class of graphs excluding at least one graph has a
low tree-depth coloring.
\end{theorem}
\begin{notation}
Following \cite{Taxi_tdepth}, we will make use of the notation
$\chi_p(G)$ for the minimum number of colors need for a vertex
coloring of $G$ such that $i<p$ parts induce a subgraph of tree-depth
at most $i$.
\end{notation}

Theorem \ref{th:ltd} relies on $p$-centered colorings, which have also been
introduced in \cite{Taxi_tdepth}:
A {\em $p$-centered coloring} of a graph $G$ is a vertex coloring
such that, for any (induced) connected subgraph $H$, either some color $c(H)$
appears exactly once in $H$, or $H$ gets at least $p$ colors.

For the sake of completeness we recall some lemmas of \cite{Taxi_tdepth}:

\begin{lemma}[\cite{Taxi_tdepth}]
\label{lem:centerlocal}
Let $G, G_0$ be graphs, let $p=\depth(G_0)$, let $c$ be a
$q$-centered coloring of $G$ where $q\geq p$.
Then any subgraph $H$ of $G$ isomorphic to $G_0$ gets at least $p$
colors in the coloring of $G$.
\qed
\end{lemma}

From this lemma follows that $p$-centered colorings induce low
tree-depth colorings:

\begin{corollary}
\label{cor:c2t}
Let $p$ be an integer, let $G$ be a graph and let $c$ be a
$p$-centered coloring of $G$.

Then $i<p$ parts induce a subgraph of tree-depth at most $i$
\end{corollary}
\begin{proof}
Let $G'$ be any subgraph of $G$ induced by $i<p$ parts. 
Assume $\depth(G')>i$. According to Lemma \ref{lem:recur}, the
deletion of one vertex decreases the tree-depth by at most one. Hence
there exists an induced subgraph $H$ of $G'$ such that
$\depth(H)=i+1\leq p$.
According to lemma \ref{lem:centerlocal} (choosing $G_0=H$), $H$ gets
at least $p$ colors, a contradiction.
\end{proof}

\begin{lemma}[\cite{Taxi_tdepth}]
\label{lem:colortw}
Let $p,k$ be integers. Then there exists an integer $N(p,k)$ such
that any graph $G$ with tree width at most $k$ has a $p$-centered coloring using
$N(p,k)$ colors.
\qed
\end{lemma}

The following lemma is proved in \cite{Taxi_tdepth} for the particular
case of proper minor closed classes of graphs and tree-width. We
shall state it here in its general form.

\begin{lemma}
\label{lem:colind}
Let $\mathcal C$ be a class of graphs. Assume that for any integer
$p\geq 1$
there exists a class of graphs $\mathcal C_p$ such that:
\begin{itemize}
\item there exists an integer $N(\mathcal C_p,p)$, such that any graph
$G\in\mathcal C_p$ has a $p$-centered coloring using at most $N(\mathcal
C_p,p)$ colors,
\item there exists an integer $C(p)$ such that any $G\in\mathcal C$
has a $C(p)$ vertex coloring such that $p$ classes induce a graph in
$\mathcal C_p$.
\end{itemize}
Then there exists an integer $X(p)$,
such that every graph in $\mathcal C$ has a $p$-centered
coloring using $X(p)$ colors.
\end{lemma}
\begin{proof}
Let $G\in\mathcal C$. According to the assumption, there exists a
vertex partition into 
$C(p)$ parts, such that any $p$ parts form a graph in $\mathcal C_p$.
This partition will be defined as a
coloring $\bar{c}: V(G) \longrightarrow \{1,2,\ldots,C(p)\}$. For
any set $P$ of $p$ parts let $G_P$ be the graph induced by all the
parts in $P$. According to the assumption, each of the
$G_P$ has $p$-centered coloring $c_P$ using $N(\mathcal C_p,p)$ colors.
Consider the following (``product'') coloring $c$ defined as
\begin{equation*}
c(v) = (\bar{c}(v), (c_P(v); |P| = p, P \subset \{1,2,\ldots,C(p)\})).
\end{equation*}

This is the product of the coloring of $G$ by $C(p)$ colors and of the
colorings of the $G_P$. This new coloring of $G$ (with $X(p)=C(p)
N(\mathcal C_p,p)^{\binom{C(p)}{p}}$ colors. Let $H$ be a connected
subgraph of $G$. Then, either $H$ gets at least $p+1$ colors, or
$V(H)$ is included in some subgraph $G_P$ of $G$ induced by $p$
parts. In the later case, some color appears exactly once in $H$.
\end{proof}

\begin{theorem}
\label{th:col}
Let $\mathcal C$ be a class of graphs having low tree-width colorings
and let $p$ be an integer. Then there exists integer $X(p)$,
such that every graph in $\mathcal C$ has a $p$-centered
coloring using $X(p)$ colors.
\end{theorem}

\begin{proof}
Let $\mathcal C_p$ be the class of graphs with tree-width at most
$(p-1)$. According to Theorem \ref{th:2tw} and
Lemma~\ref{lem:colortw}, the conditions of Lemma \ref{lem:colind} are
satisfied hence $X(p)$ exists.
\end{proof}

As a consequence we have the following equivalence of the various (seemingly unrelated) above notions:
\begin{theorem}
\label{th:equiv1}
Let $\mathcal C$ be a class of graphs. Then the following conditions
are equivalent:
\begin{itemize}
\item $\mathcal C$ has a low tree-width coloring,
\item $\mathcal C$ has a low tree-depth coloring,
\item for any integer $p$, $\{\chi_p(G): G\in\mathcal C\}$ is bounded,
\item for any integer $p$, there exists an integer $X(p)$ such that
any graph $G\in\mathcal C$ has a $p$-centered colorings using at most
$X(p)$ colors.
\end{itemize}
\end{theorem}

Our main result (Theorem \ref{thm:main}) is a non-trivial extension of this equivalence.

\section{The grad of a graph and classes with bounded expansion}
\label{sec:grad}
Recall that the {\em maximum average degree} $\mad(G)$ of a graph $G$ is
the maximum over all subgraphs $H$ of $G$ of the average degree of $H$, that is
$\mad(G)=\max_{H\subseteq G}\frac{2|E(H)|}{|V(H)|}$. The {\em
  distance} $d(x,y)$ between two vertices $x$ and $y$ of a graph is the minimum
length of a path linking $x$ and $y$, or $\infty$ if $x$ and $y$ do not
belong to the same connected component.

We introduce several notations:

\begin{itemize}
\item 
The {\em radius} $\rho(G)$ of a connected graph $G$ is:
\begin{equation*}
\rho(G)=\min_{r\in V(G)}\max_{x\in V(G)} {\rm d}(r,x)
\end{equation*}
\item
 A {\em center} of $G$ is a vertex $r$ such that $\max_{x\in V(G)}{\rm
  d}(r,x)=\rho(G)$.
\end{itemize}

\begin{definition}
Let $G$ be a graph. A {\em ball} of $G$ is a subset of vertices
inducing a connected subgraph. 
The set of all the families of balls of $G$ is noted $\mathfrak{B}(G)$.

Let $\mathcal P=\{V_1,\dotsc,V_p\}$ be a family of balls of $G$.
\begin{itemize}
\item
 The {\em radius} $\rho(\mathcal P)$ of $\mathcal P$ is
 $\rho(\mathcal P)=\max_{X\in \mathcal P}\rho(G[X])$ 
\item The {\em complexity} of $\mathcal P$
  is $\zeta(\mathcal P)=\max_{v\in V(G)}\card{\{i: v\in V_i\}}$.
\item
The {\em quotient} $G/\mathcal P$ of $G$ by 
  $\mathcal P$ is a graph with vertex set $\{1,\dotsc,p\}$ and edge
  set $E(G/\mathcal P)=\{\{i,j\}: (V_i\times V_j)\cap
  E(G)\neq\emptyset\text{ or }V_i\cap V_j\neq\emptyset\}$.
\end{itemize}
\end{definition}

We introduce several invariants that refine the notion of maximum
average degree:
\begin{definition} The {\em greatest reduced average density (grad) of $G$
      with rank $r$ and complexity $c$} is
$$\gdens{c}{r}(G)=\max_{\substack{\mathcal P\in\mathfrak{B}(G)\\
\rho(\mathcal P)\leq r, \zeta(\mathcal P)\leq
    c}}\frac{|E(G/\mathcal P)|}{|\mathcal P|}.$$
For the sake of simplicity, we also define:
\begin{itemize}
\item  The {\em grad
      of $G$ with rank $r$}:  $$\rdens{r}(G)=\gdens{1}{r}(G)=
\max_{\substack{\mathcal P\in\mathfrak{B}(G)\\
\rho(\mathcal P)\leq r,
    \zeta(\mathcal P)=1}}\frac{|E(G/\mathcal P)|}{|\mathcal P|}$$

\item   The {\em grad of $G$}: $$\dens(G)=\max_{r}\rdens{r}(G)=\max_{H\preceq
    G}\frac{|E(H)|}{|V(H)|}$$
\end{itemize}
\end{definition}

Notice that we have:
\begin{equation}
  \frac{\mad(G)}{2}=\rdens{0}(G)\leq\rdens{1}(G)\leq\dotsb\leq\rdens{\rho(G)}(G)=\dens(G)
\end{equation}
and that $\dens(G)$ is related to the Hadwiger number $h(G)$ of $G$
(that is the maximum order of a complete graph which is a minor of $G$)
by:
\begin{equation}
\frac{h(G)-1}{2}\leq\dens(G)\leq O(h(G)\sqrt{\log h(G)}),
\end{equation}
\begin{proof}
Let $h=h(G)$. As $K_h$ is a $(h-1)$-regular minor of $G$ ,
$\frac{h-1}{2}\leq\dens(G)$. Moreover, there exists a constant $C$
such that if $\dens(G) > C (h+1)\sqrt{\log (h+1)}$ then $G$ has a minor with
minimum degree at least $\gamma (h+1)\sqrt{\log (h+1)}$ hence a minor
$K_{h+1}$ as proved by  Kostochka \cite{Ko} and Thomason \cite{TH1}
(extending earlier work of Mader \cite{M}; see \cite{Th2} for an tight
value of constant $\gamma$).
\end{proof}

Also notice the following well known facts (usually expressed by means
of the maximum average degree):
\begin{fact}
\label{fact:md}
Let $G$ be a graph. Then $G$ has an orientation such that the maximum
indegree of $G$ is at most $k$ if and only if $k\geq \rdens{0}(G)$.
\end{fact}
\begin{fact}
Let $G$ be a graph. Then $G$ is $\lfloor
2\rdens{0}(G)\rfloor$-degenerated, hence 
$\lfloor 2\rdens{0}(G)+1\rfloor$-colorable.
\end{fact}

The grad actually appears to be related to low tree-depth colorings:

\begin{lemma}
\label{lem:grad_chi}
For any graph $G$ and any integer $r$:
\begin{equation}
 \rdens{r}(G)\leq (2r+1)\binom{\chi_{2r+2}(G)}{2r+2}
\end{equation}
\end{lemma}
\begin{proof}
 Consider a vertex coloring $c$ of $G$ with $N=\chi_{2r+2}(G)$ colors such that any
 $i\leq 2r+2$ colors induce a subgraph of tree-depth at most $i$. 
For any $J\in\binom{[N]}{2r+2}$, let $G_J=G[c^{-1}(J)]$ and let
$Y_J$ be a rooted forest of height $\depth(G_J)\leq 2r+2$ such that
$G_J\subseteq\clos(Y_J)$. 

Let
 $\mathcal P=\{X_1,\dotsc,X_p\}$ be a family of balls of $G$ with radius $r$ and
 complexity $1$ achieving the bound $\rdens{r}(G)$ (that is: such that
 $\rdens{r}(G)=\frac{\card{E(G/\mathcal P)}}{\card{\mathcal P}}$). 
Let
 $x_1,\dotsc,x_k$ be centers of $X_1,\dots,X_k$. 
 If $X_i$ and $X_j$ are adjacent in $G/\mathcal P$
 then there exists a path $P_{i,j}$ of length at most $2r+1$ linking $x_i$ and
 $x_j$. Let $I_{i,j}\in\binom{[N]}{2r+2}$ be such that
 $I_{i,j}\supseteq c(V(P_{i,j}))$. Then the path $P_{i,j}$ is included
 in some connected 
 component of $G_{I_{i,j}}$. It follows that there exists in $P_{i,j}$
 a vertex $v_{i,j}$ which is minimum with respect to the partial order
 defined by $Y_{I_{i,j}}$. As $\{x_i,x_j\}\subseteq
 V(P_{i,j})\subseteq X_i\cup X_j$ and as $X_i\cap X_j=\emptyset$
 (because $\zeta(\mathcal P)=1$), $v_{i,j}$ either belongs to $X_i$ or
 to $X_j$. Depending on the case, $v_{i,j}$ is a vertex of $X_i$ which
 is an ancestor of $x_j$ in $Y_{I_{i,j}}\cap X_i$ or a vertex of $X_j$ which
 is an ancestor of $x_i$ in $Y_{I_{i,j}}\cap X_j$.
 Thus:
\begin{align*}
p\rdens{r}(G)&\leq
\sum_{I\in\binom{[N]}{2r+2}}
    \sum_{1\leq i\leq p}\sum_{\substack{1\leq j\leq
    p\\ j\neq i}} \quad\card{\{v:
  v\text{ ancestor of }x_i\text{ in }Y_I\cap X_j\}}\\
&\leq
\sum_{I\in\binom{[N]}{2r+2}}\sum_{i=1}^{p} \quad\card{\{v:
  v\text{ ancestor of }x_i\text{ in }Y_I\}}\\
&\leq \binom{N}{2r+2}\times p\times(2r+1)
\intertext{Hence}
\rdens{r}(G)&\leq(2r+1)\binom{N}{2r+1}
\end{align*}
 \end{proof}

This lemma motivates  the following definition:

\begin{definition}
A class of graphs $\mathcal C$ has {\em bounded expansion} if there
exists a function $f:\bbbn\rightarrow\bbbr$ such that for every graph $G \in \mathcal C$ and every $r$ holds
\begin{equation}
 \rdens{r}(G)\leq f(r)
\end{equation}
\end{definition}

\begin{theorem}
\label{th:tw2be}
If a class $\mathcal C$ has low tree-width colorings then
$\mathcal C$ has bounded expansion.
\end{theorem}
\begin{proof}
As low tree-width colorings and low tree-depth colorings are
equivalent, the theorem is a direct consequence of Lemma \ref{lem:grad_chi}.
\end{proof}

The main theorem of this paper may be seen as a converse of Theorem \ref{th:tw2be}.

\section{Grad stability over lexicographic product}
\label{sec:lex}
Let $G,H$ be graphs. The {\em lexicographic product} $G\bullet H$ is
defined by $V(G\bullet H)=V(G)\times V(H)$ and $E(G\bullet
H)=\{\{(x,y),(x',y'): \{x,y\}\in E(G)\text{ or }x=x'\text{ and
}\{y,y'\}\in E(H)\}$.

Let us note at this place that the lexicographic product 
(or blowing up of vertices) is an operation which is incompatible 
with the minors.  One can see easily that every graph is a minor of a 
graph of the form $G\bullet K_2$ for a planar graph $G$.
But the lexicographic product is naturally related to the notion of
complexity we have introduced for grad:

\begin{lemma}
\label{lem:lex}
  For any graph $G$ and any integers $c,r$, we have:
  \begin{equation*}
    \gdens{c}{r}(G)=\rdens{r}(G\bullet K_c)
  \end{equation*}
\end{lemma}
\begin{proof}
  Let $\mathcal P=\{V_1,\dots,V_p\}$ be a ball family of $G$
  with complexity $c=\zeta(\mathcal P)$ and radius $r=\rho(\mathcal
  P)$.  As $\zeta(\mathcal P)=c$ there exists a function
  $f:V(G)\times\{1,\dots,p\}\rightarrow\{1,\dots,c\}$ such that if
  $x\in V_i\cap V_j$ then $f(x,i)\neq f(x,j)$.

  For $1\leq i'\leq p$, define $V_i'=\{(x,f(x,i)): x\in V_i\}$.  
Then $\mathcal
  P'=\{V_1',\dots,V_p'\}$ has radius $r$ and complexity $1$.
  Moreover, $G/\mathcal P$ is obviously isomorphic to a subgraph of
  $(G\bullet K_c)/\mathcal P'$. It follows that $\rdens{r}(G\bullet
  K_c)\geq \gdens{c}{r}(G)$.

  Conversely, let $\mathcal P'=\{V_1',\dots,V_q'\}$ be a ball family
  of $G\bullet K_c$, define the ball family $\mathcal
  P=\{V_1,\dots,V_q\}$ of $G$ by $x\in V_i$ if there exists
  $\alpha\in\{1,\dots,c\}$ such that $(x,\alpha)\in V_i'$. Then
  $\rho(\mathcal P)\leq \rho(\mathcal P')$ and $\zeta(\mathcal P)\leq
  c$. It follows that $\gdens{c}{r}(G)\geq\rdens{r}(G\bullet K_c)$.
\end{proof}

The remaining of the section will be dedicated to the proof of the
following key lemma:
\begin{lemma}
\label{lem:gdens}
There exist polynomials $P_i\ (i\geq 0)$ such that
for any graph $G$ and integers
$r$ and $c$:
\begin{equation}
\gdens{c}{r}(G)\leq P_r(c,\rdens{r}(G))
\end{equation}
\end{lemma}

In the following, a directed graph $\vG$ may not have a loop and for
any two of its vertices $x$ and $y$, $\vG$ includes at most one arc
from $x$ to $y$ and at most one arc from $y$ to $x$.

If a directed path $\vP$ has starting vertex $x$ and end vertex $y$,
we note $\dpath{x}{\vP}{y}$.

If $\dpath{x}{\vP_1}{z}$, $\dpath{y}{\vP_2}{z}$ and if 
no internal vertex or edges of $\vP_1$ belongs to $\vP_2$ nor the
converse, we note $\lambd{x}{\vP_1}{z}{\vP_2}{y}$. In such a case,
either $\vP_1\cup\vP_2$ is a path, or  $\vP_1\cup\vP_2$ is a cycle and
$x=y$. Moreover, if $x\neq y$, $\card{\vP_1}\leq a$ and
$\card{\vP_2}\leq b$, we say that $y$ is {\em $(a,b)$-reachable} from
$x$. 
\begin{definition}
Let $\vG$ be a directed graph, let $a,b$ be integers.
A set $\vL$ of arcs with endpoints in $V(\vG)$ is an 
{\em $(a,b)$-augmentation} of $\vG$ if, for any $x,y\in V(\vG)$
with $y$ $(a,b)$-reachable from $x$, either $(x,y)$ or $(y,x)$ belongs
to $\vL$.

The {\em maximum indegree} of $\vL$ is
\begin{equation*}
  \md(\vL)=\max_{y\in V(\vG)}\card{\{x\in V(\vG):
    (x,y)\in\vL\}}
\end{equation*}
\end{definition}

Notice that if $a$ or $b$ is at least $1$, $E(\vG)$ is obviously included in any
$(a,b)$-augmentation of $\vG$.

\begin{lemma}
  \label{lem:vcol}
  Let $\vG$ be a directed graph, let $a,b$ be integers and let
  $\vL$ be an $(a,b)$-augmentation of $\vG$.

  Then there exists a vertex coloring $\vcol{\vL}$ using at most
  $2\md(\vL)+1$ colors such that for any vertex $x$,
  $\vcol{\vL}(y)\neq\vcol{\vL}(x)$ for any vertex $y$ $(a,b)$-reachable from
  $x$.
\end{lemma}
\begin{proof}
  Let $\vH$ be the directed graph with vertex set $\vG$ and arc set
  $\vL$. If $y$ is $(a,b)$-reachable from $x$ in $\vG$ then $(x,y)$
  or $(y,x)$ belongs to $E(\vH)$. As $\vH$ has maximum indegree
  $\md(\vL)$, it is $(2\md(\vL)+1)$-choosable. Any proper
  coloration of $\vH$ will do.
\end{proof}
\begin{lemma}
  \label{lem:ecol}
  Let $\vG$ be a directed graph with maximum indegree $\md(\vG)$, let
  $a,b$ be integers and let $\vL$ be an $(a,b)$-augmentation of $\vG$.

 Then there exists
  an edge coloring $\ecol{\vL}$ using at most $(2\md(\vL)+1)\md(\vG)$
  colors such that for any
$\lambd{x}{\vP_1}{z}{\vP_2}{y}$
 with $\card{\vP_1}\leq a+1$ and $\card{\vP_2}\leq b+1$, 
all the edges of $\vP_1\cup\vP_2$ get different colors.
\end{lemma}
\begin{proof}
  Consider an edge coloring $c_0$ such that two edges having the
  same end vertex have different colors (this is achieved with
  $\md(\vG)$ colors) and the vertex coloring $\vcol{\vL}$ defined in
  Lemma \ref{lem:vcol}.  Then for any arc $e=(x,y)$ define
  $\ecol{\vL}(e)=(c_0(e),\vcol{\vL}(y))$. Then if $e=(x,y)$ and $f=(x',y')$
  are two different arcs in $\vP_1\cup\vP_2$ where
   either $y\neq y'$ thus $y'$ is
  $(a,b)$-reachable from $y$ or $y$ is $(a,b)$-reachable from $y'$
  hence $\vcol{\vL}(y')\neq\vcol{\vL}(y)$, or $y=y'$ hence
  $c_0(e)\neq c_0(f)$.
\end{proof}
\begin{notation}
Let $\Upsilon$ be an arc-coloring of a directed graph $\vG$ and let
$\vP$ be a directed path of $\vG$ of length $l$. 
We note $\Upsilon(\vP)=\vec{\alpha}=(\alpha_1,\dotsc,\alpha_l)$ the
sequence of the colors $\Upsilon(e)$ of the arcs of $\vP$, taken in the
order in which they appear on $\vP$.
\end{notation}

\begin{lemma}
  \label{lem:tree0}
  Let $\vG$ be a directed graph with maximum indegree $\md(\vG)$, let
  $a,b$ be integers and let $\vL$ be an $(a,b)$-augmentation of
  $\vG$. Let $\ecol{\vL}$ be the edge coloring defined in
  Lemma~\ref{lem:ecol}. 

  Let $\vP_1,\vP_2$ be two directed paths of length
  $l\leq\max(a,b)+1$, such that the initial vertex of one of them is
  different from the end vertex of the other one.
  If $\ecol{\vL}(\vP_1)=\ecol{\vL}(\vP_2)$ then 
  either $\vP_1$ and $\vP_2$ do
  not intersect, or they share the same initial vertex and 
  there exists $0\leq a\leq l$ such that $\vP_1$ and
  $\vP_2$ share their $a$ first edges and do not intersect
  thereafter.
\end{lemma}
\begin{proof}
Without loss of generality, we may assume $a\geq b$. 
Let $\vec{\alpha}=\ecol{\vL}(\vP_1)$.
Assume there exists a vertex $v$ having one incoming edge in $\vP_1$
(the $i$th of $\vP_1$, hence colored $\alpha_i$) and one
(different) incoming edge in $\vP_2$ (the $j$th of $\vP_2$, hence
colored $\alpha_j$). Without loss of generality, we may assume $i\geq
j$. 
Then the $(j+1)$th vertex $u$ of $\vP_1$ has an incoming edge in
$\vP_1$ colored $\alpha_j$ and belong to the initial subpath of
$\vP_1$ ending at $v$. It follows that $v$ is $(a,0)$ reachable from
$u$. Hence an incoming edge of $u$ may not have the same color of an
incoming edge of $v$, contradiction.
Similarly, the initial vertex of one of the path may not be internal
to the second one.
As the case where the initial vertex of one of the path is the end
vertex of the other one, we conclude that either the two paths do not
intersect or they share their $a$ first edges.
\end{proof}
\begin{lemma}
  \label{lem:tree}
  Let $\vG$ be a directed graph with maximum indegree $\md(\vG)$, let
  $a,b$ be integers and let $\vL$ be an $(a,b)$-augmentation of
  $\vG$. Let $\ecol{\vL}$ be the edge coloring defined in
  Lemma~\ref{lem:ecol}. Let $\vec{\alpha}$ be a sequence of 
  $l\leq \max(a,b)+1$ distinct edge colors. Then the union
  $T_{\vL}(\vec{\alpha})$ of all the directed paths $\vP$ such that
  $\ecol{\vL}(\vP)=\vec{\alpha}$ is a directed rooted forest.
\end{lemma}
\begin{proof}
This is a direct consequence of Lemma~\ref{lem:tree0}.
\end{proof}
\begin{lemma}
  \label{lem:nocross}
  Let $\vG$ be a directed graph with maximum indegree $\md(\vG)$, let
  $a\geq b$ be integers and let $\vL$ be an $(a,b)$-augmentation of
  $\vG$. Let $\ecol{\vL}$ be the edge coloring defined in
  Lemma~\ref{lem:ecol}. Let $\vec{\alpha}$ and $\vec{\beta}$ be
  sequences of respective lengths $p\leq a+1$ and $q\leq b+1$. 
  Let $\Pi_{\vL}(\vec{\alpha},\vec{\beta})$ be the union of all the
  $\vP_1\cup\vP_2$ where $\ecol{\vL}(\vP_1)=\vec{\alpha}$,
  $\ecol{\vL}(\vP_2)=\vec{\beta}$ and there exists three distinct
  vertices $x,y,z$ so that $\lambd{x}{\vP_1}{z}{\vP_2}{y}$.

  Then a directed tree $Y_1$ in
  $\Pi_{\vL}(\vec{\alpha},\vec{\beta})\cap T_{\vL}(\vec{\alpha})$ and
  a directed tree $Y_2$ in $\Pi_{\vL}(\vec{\alpha},\vec{\beta})\cap
  T_{\vL}(\vec{\beta})$ with different roots may only intersect at a
  leaf of both of them.
\end{lemma}
\begin{proof}
Let $r_1,r_2$ be the roots of $Y_1$ and $Y_2$.
If $Y_1$ and $Y_2$ intersects, there exists 
$\lambd{r_1}{\vP_1}{z}{\vP_2}{y}$ and 
$\lambd{x'}{\vP_1'}{z'}{\vP_2'}{r_2}$, so that
$\ecol{\vL}(\vP_1)=\ecol{\vL}(\vP_1')=\vec{\alpha}$,
$\ecol{\vL}(\vP_2)=\ecol{\vL}(\vP_2')=\vec{\beta}$,
and $\vP_2'$ intersects $\vP_1$ at a vertex $v$ (up to an exchange of
 $Y_1$ and $Y_2$). As $r_1\neq r_2$, $v$ has in $\vP_2$ an incoming
 edge $e$ of color $\beta_i$ for some $1\leq i\leq b+1$. Let
 $w$ be the vertex of $\vP_2$ having in $\vP_2$ an incoming edge of
 color $\beta_i$. If $w\neq v$, we are led to a contradiction,
 according to Lemma~\ref{lem:ecol}, as $w$
 is then $(p,q)$-reachable from $v$. Hence $v=w$ and $v$ is the end
 vertex of $\vP_1$ and $\vP_2$. Thus $v$ is also the end vertex of
 $\vP_1'$ and $\vP_2'$. It follows that $v$ is a leaf of both $Y_1$
 and $Y_2$.
\end{proof}
\begin{lemma}
\label{lem:augm}
  Let $\vG$ be a directed graph with maximum indegree $\md(\vG)$, let
  $r$ be an integer and let $\vL$ be an $(r,r-1)$-augmentation of
  $\vG$.

  Then $\vL$ may be extended into an $(r+1,r)$-augmentation $\vL'$
  such that
  $\md(\vL')\leq \md(\vL)+((2\md(\vL)+1)\md(\vG))^{2r+1}\rdens{r}(G)$.
\end{lemma}
\begin{proof}
Let $\ecol{\vL}$ be the edge coloring defined in
  Lemma~\ref{lem:ecol}. 

  For two sequences $\vec{\alpha}$ and $\vec{\beta}$
  of respective lengths $p\leq r+1$ and $q\leq r$,
  let $\Pi_{\vL}(\vec{\alpha},\vec{\beta})$ be the union of all the
  $\vP_1\cup\vP_2$ where $\ecol{\vL}(\vP_1)=\vec{\alpha}$,
  $\ecol{\vL}(\vP_2)=\vec{\beta}$ and there exists three distinct
  vertices $x,y,z$ so that $\lambd{x}{\vP_1}{z}{\vP_2}{y}$. Also,
  let $G_{\vec{\alpha},\vec{\beta}}$ be the graph obtained from $G$ by
  contracting all the edges of $\Pi_{\vL}(\vec{\alpha},\vec{\beta})$
  but those colored $\alpha_p$.

  Let $x,y$ be vertices of $G$ so that $y$ is
  $(r+1,r)$-reachable from $x$, as witnessed by
  $\lambd{x}{\vP_1}{z}{\vP_2}{y}$. Let 
  $\vec{\alpha}=\ecol{\vL}(\vP_1)$ and
  $\vec{\beta}=\ecol{\vL}(\vP_2)$.
  The vertices $x,y$ are the roots of directed trees in 
  $\Pi_{\vL}(\vec{\alpha},\vec{\beta})\cap T_{\vL}(\vec{\alpha})$ and
  $\Pi_{\vL}(\vec{\alpha},\vec{\beta})\cap T_{\vL}(\vec{\beta})$,
  respectively, hence to two adjacent distinct vertices in
  $G_{\vec{\alpha},\vec{\beta}}$.
  Similarly, two distinct vertices of $G_{\vec{\alpha},\vec{\beta}}$
  adjacent by an edge of color $\alpha_p$ (where
  $p=\card{\vec{\alpha}}$) correspond uniquely to the roots of a
  tree in $\Pi_{\vL}(\vec{\alpha},\vec{\beta})\cap T_{\vL}(\vec{\alpha})$ and
  $\Pi_{\vL}(\vec{\alpha},\vec{\beta})\cap T_{\vL}(\vec{\beta})$,
  respectively.

  It follows that there exists an $(r+1,r)$-augmentation $\vL'$ of $\vG$
  extending $\vL$ such that  $$\md(\vL')-\md(\vL)\leq \sum_{\substack{\card{\vec{\alpha}}\leq r+1\\
      \card{\vec{\beta}}\leq r}}
  \rdens{0}(G_{\vec{\alpha},\vec{\beta}})\leq
 ((2\md(\vL)+1)\md(\vG))^{2r+1}\rdens{r}(G)
$$

\end{proof}
\begin{lemma}
\label{lem:augm2}
For any integer $r$, there exists a polynomial $\Phi_r$ such that any
directed graph $\vG$ has a $(r+1,r)$-augmentation $\vL$, where
$\md(\vL)\leq \Phi_r(\md(\vG),\rdens{r}(G))$, where $G$ is the
underlying simple graph of $\vG$.
\end{lemma}
\begin{proof}
This is a direct consequence of Lemma~\ref{lem:augm}.
\end{proof}
\begin{proof}[\bf Proof of Lemma \ref{lem:gdens}]
Define $P_r(x,y)=\Phi_r(x+y,y)$.

Consider a family $\mathcal P$ of balls of $G$ with radius at most
$r$ and complexity at most $c$. We construct a directed graph $\vG$
with underlying undirected graph $G$. Recall that $\vG$ may have, for
each edge of $G$, one arc in each direction.
First we orient the edges of $G$ with
indegree $\rdens{0}(G)$ (thus obtaining one arc per edge).
 For 
each $X\in\mathcal P$, let $v$ be the center of $G[X]$. Let $Y$ be a
minimum distance tree of $G[X]$ with root $v$. If $\vG$ does not
include the arcs corresponding to an orientation of $Y$ from its root
$v$, we add the missing arcs. We also add if necessary
all the arcs going from a leaf of $Y$ to a vertex out of
$X$.

Notice that the vertices of $\vG$
have indegree at most $\rdens{0}(G)+c$. Moreover, if $r_1,r_2$
are the roots of the trees $Y_1$ and $Y_2$ corresponding to
some parts $X_1,X_2\in\mathcal P$ which are adjacent in $G/\mathcal P$ then
$r_2$ is $(r+1,r)$-reachable from $r_1$ in $\vG$ (by a directed path of
length at most $r$ in $Y_1$, possibly followed by an arc between the
parts and a directed path of length at most $r$ in $Y_2$ in opposite
direction). Hence $r_1$ and $r_2$ are adjacent in any
$(r+1,r)$-augmentation of $\vG$. According to Lemma~\ref{lem:augm2},
there exists such an augmentation $\vL$ with
$\md(\vL)\leq \Phi_r(\rdens{0}(G)+c,\rdens{r}(G))$.
As $G/\mathcal P$ is isomorphic to a subgraph of the graph with vertex
set $V(G)$ and edge set $\vL$. As this subgraph has an orientation
with indegree at most $\md(\vL)$ we have, according to Fact
\ref{fact:md} and Lemma \ref{lem:augm2}:
\begin{equation*}
\gdens{c}{r}(G)=\rdens{0}(G/\mathcal P)\leq \md(\vL)\leq
\Phi_r(\rdens{0}(G)+c,\rdens{r}(G))\leq P_r(c,\rdens{r}(G)).
\end{equation*}

\end{proof}
\section{Transitive fraternal augmentation}
\label{sec:aug}
\begin{definition}
Let $\vG$ be a directed graph. A {\em $1$-transitive fraternal augmentation} of $\vG$
is a directed graph $\vH$ with the same vertex set, including all the
arcs of $\vG$ and such that, for any vertices $x,y,z$,
\begin{itemize}
\item if $(x,z)$ and $(z,y)$ are arcs of $\vG$ then $(x,y)$ is an arc
  of $\vH$ ({\em transitivity}),
\item  if $(x,z)$ and $(y,z)$ are arcs of $\vG$ then $(x,y)$ or $(y,x)$ is an arc
  of $\vH$ ({\em fraternity}).
\end{itemize}

A {\em transitive fraternal augmentation} of a directed graph $\vG$ is a sequence
$\vG=\vG_1\subseteq\vG_2\subseteq\dotsb\subseteq \vG_i\subseteq\vG_{i+1}\subseteq\dotsb$,
such that $\vG_{i+1}$ is a $1$-transitive fraternal augmentation of $\vG_i$ for any
$i\geq 1$.
\end{definition}

The main key lemma here is that the notion of classes of bounded
 expansion is stable under
 $1$-fraternal augmentations.
More precisely:

\begin{lemma}
\label{lem:faug}
Let $\vG$ be a directed graph and let $\vH$ be a $1$-transitive fraternal augmentation
of $\vG$. Then
\begin{equation}
\gdens{c}{r}(H)\leq \gdens{c(\md(\vG)+1)}{2r+1}(G)\leq P_{2r+1}(c(\md(\vG)+1),\rdens{2r+1}(G)).
\end{equation}
\end{lemma}
\begin{proof}
Consider a ball family $\mathcal P=\{V_1,\dotsc,V_p\}$ of $H$ with radius at most $r$ and
complexity $c$. 
Let $\mathcal P'=\{V_1',\dotsc,V_p'\}$, where
 $V_i'=V_i\cup\{z: \exists x\in V_i, (x,z)\in E(\vG)\}$.
Then for any $x,y\in V_i$ which are adjacent in $H$, either $x$ and
$y$ are adjacent in $G$ or there exists $z\in V_i'$ so that $\{x,z\}$
and $\{y,z\}$ are edges of $G$. Hence $V_i'$ is a ball of $G$ with
radius at most $2r+1$. Any vertex $v$ of $G$ belongs to a most
$c+\md{\vG}$ balls of $\mathcal P'$ for $v$ belongs to $V_i'$ if and
only if either $v$ belongs to $V_i$ (there are at most $c$ such $V_i$)
or there exists an arc from a vertex $z\in V_i$ to $v$ in $\vG$ (there
are at most $\md(\vG)$ such $z$ hence at most $c\md(\vG)$ such $V_i$).
Hence the complexity of $\mathcal P'$ is at most $c(\md(\vG)+1)$.
As $H/\mathcal P$ is isomorphic to a subgraph of $G/\mathcal P'$
$\card{E(H/\mathcal P)}\leq\card{E(G/\mathcal P')}$ thus
$\gdens{c}{r}(H)=\frac{\card{E(H/\mathcal P)}}{\card{\mathcal P}}\leq
\frac{\card{E(G/\mathcal P')}}{\card{\mathcal P'}}\leq\gdens{c(\md(\vG)+1)}{2r+1}(G)$.
We conclude using Lemma~\ref{lem:gdens}.
\end{proof}
\begin{corollary}
\label{cor:faugm}
There exists polynomials $Q_i\ (i\geq 1)$, such that
any directed graph $\vG$ has a transitive fraternal augmentation
$\vG=\vG_1\subseteq\vG_2\subseteq\dotsb\subseteq \vG_i\subseteq\dotsb$
where 
\begin{equation}
\md(\vG_i)\leq Q_i(\md(\vG),\rdens{2^{i+1}-1}(G))
\end{equation}
\end{corollary}
We also deduce:

\begin{corollary}
Let $\mathcal C$ be a class with bounded expansion. Then there exists
a function $g$ such that each graph $G\in\mathcal C$ has a transitive
fraternal augmentation $\vG=\vG_1\subseteq\vG_2\subseteq\dotsb\subseteq \vG_i\subseteq\dotsb$
where $\md(\vG_i)\leq g(i)$.
\end{corollary}

\section{Back to $p$-centered colorings}
\label{sec:back}
The aim of this section is to prove that transitive fraternal
augmentations allow us to construct $p$-centered colorings.

\begin{lemma}
\label{lem:pt}
Let $N(p,t)=1+(t-1)(2+\lceil\log_2 p\rceil)$,
let $\vG$ be a directed graph and let
$\vG=\vG_1\subseteq\vG_2\subseteq\dotsb\subseteq \vG_i\subseteq\dotsb$
be a transitive fraternal augmentation of $\vG$.

Then $\vG_{N(p,\depth(G))}$ either includes an acyclically oriented clique of
size $p$ or a rooted directed tree $\vY$ such that $G\subseteq\clos(Y)$
and $\clos(\vY)\subseteq\vG_{N(p,\depth(G))}$.
\end{lemma}
\begin{proof}
We fix the integer $p$ and prove the lemma by induction on $t=\depth(\vG)$.
The base case $t=1$ corresponds to a graph without edges, for which
the property is obvious. Assume the lemma has been proved for directed
graphs with tree-depth at most $t$ and let $\vG$ be a directed graph
with tree-depth $t+1$. As we may consider each connected component of
$\vG$ independently, we may assume that $\vG$ is connected. Then there
exists a vertex $s\in V(\vG)$ such that the connected components
$\vH_1,\dotsc,\vH_k$ of $G-s$ have tree-depth at most $t$. As
$\vH_i=\vG_1[V(\vH_i)]\subseteq\dotsb\subseteq\vG_j[V(\vH_i)]\subseteq\dotsb$ is a transitive fraternal
augmentation of $\vH_i$ we have, according to the induction
hypothesis, that, for each $1\leq i\leq k$, there exists in $\vH_i$
either an acyclically oriented clique of size $p$ or a rooted tree
$\vY_i$ rooted at $r_i$ such that $H_i\subseteq\clos(Y_i)$
and $\clos(\vY_i)\subseteq\vG_{N(p,\depth(G))}[V(\vH_i)]$. If the
first case occurs for some $i$, then $\vG$ includes an acyclically
oriented clique of size $p$. Hence assume it does not. As $\vG$ is
connected, the vertex $s$ has at least a neighbor $x_i$ in $\vH_i$
(for each $1\leq i\leq k$). Let $x$ be any neighbor of $s$ in
$\vH_i$. If $y$ is an ancestor of $x$ in $\vY_i$, 
$(y,x)$ is an arc of $\vG_{N(p,t)}$ hence $s$ and $y$ are
adjacent in $\vG_{N(p,t)+1}$. Moreover, if $(x,s)$ is an arc of
$\vG_{N(p,t)}$ then $(y,s)$ is an arc of $\vG_{N(p,t)+1}$. Let $D_i$
be the subset of $V(\vH_i)$ of the vertices $x$ such that $(x,s)$
belongs to $\vG_{N(p,t)}$ and of their ancestors in $\vY_i$ and let
$D=\bigcup_{i=1}^k D_i$. Then $D$ includes a clique in
$\vG_{N(p,t)+2}$. Thus there exists a
directed Hamiltonian path $\vP$ in $\vG_{N(p,t)+2}[D]$.

Let $r$ be the start vertex of $\vP$.
Define $\pi: V(G)-r\rightarrow V(G)$ as follows:
\begin{itemize}
\item if $x\in D$, the $\pi(x)$ is the predecessor $y$ of
  $x$ in $\vP$ (the arc $(y,x)$ belongs to $\vG_{N(p,t)+2}$);
\item otherwise, if $x=s$, $\pi(x)$ is the end vertex $y$ of
  $\vP$ (the arc $(y,x)$ belongs to $\vG_{N(p,t)+1}$);
\item otherwise, if $x=r_i$ then  $\pi(x)=s$
(the arc $(s,r_i)$ belongs to $\vG_{N(p,t)+2}$);
\item otherwise, if the father of $x\in V(\vH_i)\setminus D$ does not
  belong to $D$, then $\pi(x)$ is the father of $x$
  in $\vY_i$;
\item otherwise, if no descendant of $x$ in $\vY_i$ has an arc coming
  from $s$ in $\vG_{N(p,t)+1}$, $\pi(x)$ is the
  father of $x$ in $\vY_i$;
\item otherwise, $\pi(x)=s$ (the arc $(s,x)$
  belongs to $\vG_{N(p,t)+2}$).
\end{itemize}
It is easily checked that the so defined ``father mapping'' $\pi$ actually
defines a directed rooted tree $\vY$ of  $\vG_{N(p,t)+2}$ with root
$r$ and that $G\subseteq\clos(\vY)$. Moreover, either $\vY$ has height
at least $p$ and  $\vG_{N(p,t)+2+\lceil\log_2 p\rceil}$ includes an acyclically
oriented clique of size $p$ or $\clos(\vY)\subseteq
\vG_{N(p,t)+2+\lceil\log_2 p\rceil}$.
As $N(p,t+1)=N(p,t)+2+\lceil\log_2 p\rceil$, the induction follows.
\end{proof}
\begin{lemma}
\label{lem:p}
Let $p$ be an integer, let $\vG$ be a directed graph and let
$\vG=\vG_1\subseteq\vG_2\subseteq\dotsb\subseteq \vG_i\subseteq\dotsb$
be a transitive fraternal augmentation of $\vG$.
Then either $\vG_{N(p,p)}$ includes an acyclically oriented clique of
size $p$ or $\depth(G)\leq p-1$ and there
exists in $\vG_{N(p,p)}$ a rooted directed tree $Y$ so that
$G\subseteq\clos(Y)$ and $\clos(\vY)\subseteq \vG_{N(p,p)}$.
\end{lemma}
\begin{proof}
If $\depth(G)>p$ we may consider a connected subgraph of $H$ of tree-depth
$p$. According to Lemma~\ref{lem:pt}, there will exists in
$\vG_{N(p,p)}[V(H)]$ an acyclically
oriented clique of size $p$ or a rooted directed tree $\vY$ so that
$H\subseteq\clos(Y)$ and $\clos(\vY)\subseteq \vG_{N(p,p)}[V(H)]$. In
the later case, if $\depth(G)=p$ then the height of $\vY$ is at least
$\depth(H)=p$ and $\clos(\vY)$ includes an acyclically oriented clique
of size $p$.
\end{proof}

\begin{corollary}
\label{cor:aug_chi}
Let $R(p)=1+(p-1)(2+\lceil\log_2 p\rceil)=O(p\log_2 p)$.

For any graph $G$, for any transitive fraternal augmentation
$\vG=\vG_1\subseteq\vG_2\subseteq\dotsb\subseteq \vG_i\subseteq\dotsb$
of $G$ and for any integer $p$:
\begin{equation}
 \chi_p(G)\leq 2\md(\vG_{R(p)})+1
\end{equation}
\end{corollary}
And also:
\begin{corollary}
\label{cor:aug_chi2}
Let $\mathcal C$ be a class of graphs.
Assume there exists a function $f$ such that each graph $G\in\mathcal
C$ has a transitive fraternal augmentation
$\vG=\vG_1\subseteq\vG_2\subseteq\dotsb\subseteq \vG_i\subseteq\dotsb$
such that $\md(\vG_i)\leq f(i)$. Then, for any integer $p$ there
exists an integer $X(p)$ such that every $G\in\mathcal C$ has a
$p$-centered coloring using at most $X(p)$ colors.
\end{corollary}
\section{Conclusion}
\label{sec:cncl}
All previous results are gathered in the
following equivalence:
\begin{theorem}
\label{thm:main}
Let $\mathcal C$ be a class of graphs. The following conditions are
equivalent:
\begin{itemize}
\item $\mathcal C$ has low tree-width colorings,
\item $\mathcal C$ has low tree-depth colorings,
\item for any integer $p$, $\{\chi_p(G): G\in\mathcal C\}$ is bounded,
\item for any integer $p$, there exists an integer $X(p)$ such that
any graph $G\in\mathcal C$ has a $p$-centered colorings using at most
$X(p)$ colors,
\item $\mathcal C$ has bounded expansion,
\item for any integer $c$, the class  $\mathcal C\bullet
  K_c=\{G\bullet K_c: G\in\mathcal C\}$ has bounded expansion,
\item for any integer $k$, the class $\mathcal C'$ of the
$1$-transitive fraternal augmentations of directed graphs $\vG$ with
$\md(\vG)\leq k$ and $G\in\mathcal C$ form a class with bounded expansion,
\item there exists a function $F$ such that
any orientation $\vG$ of a graph $G\in\mathcal C$ has a
  transitive fraternal augmentation
  $\vG=\vG_1\subseteq\vG_2\subseteq\dotsb\subseteq\vG_i\subseteq\dotsb$
  where $\md(\vG_i)\leq F(\md(\vG),i)$,
\item there exists a function $f$ such that
  any graph $G\in\mathcal C$ has a
  transitive fraternal augmentation
  $\vG=\vG_1\subseteq\vG_2\subseteq\dotsb\subseteq\vG_i\subseteq\dotsb$
  where $\md(\vG_i)\leq f(i)$.
\end{itemize}
\end{theorem}

Now that we know that  bounded expansion is the more
general condition for low tree-depth coloring to exist and that low
tree-width coloring (although seemingly weaker) does not relax this
condition, we may wonder what may be the weakest coloring condition
equivalent to low tree-width coloring. It appears that this is a
direct consequence of Lemma~\ref{lem:colind} and
Theorem~\ref{thm:main}:

\begin{corollary}
Let $\mathcal C$ be a class of graphs. Then $\mathcal C$ has bounded
expansion if, and only if, for every integer $p\geq 1$ there exists a
class of graphs $\mathcal C_p$ and an integer $C(p)$ such that:
\begin{itemize}
\item $\mathcal C_p$ has bounded expansion,
\item any graph $G\in\mathcal C$ has a $C(p)$ vertex-coloring such
that any $p$ parts induce a graph in $\mathcal C_p$.
\end{itemize}
\end{corollary}


\providecommand{\bysame}{\leavevmode\hbox to3em{\hrulefill}\thinspace}
\providecommand{\MR}{\relax\ifhmode\unskip\space\fi MR }
\providecommand{\MRhref}[2]{%
  \href{http://www.ams.org/mathscinet-getitem?mr=#1}{#2}
}
\providecommand{\href}[2]{#2}

\end{document}